\documentclass[12pt]{article} 

\usepackage[utf8]{inputenc} 

\usepackage{geometry} 
\usepackage{graphicx} 

\usepackage{booktabs} 
\usepackage{array} 
\usepackage{paralist} 
\usepackage{verbatim} 
\usepackage{subfig} 

\usepackage{amsthm}
\usepackage{amsmath}
\usepackage{amssymb}
\usepackage{graphicx}
\usepackage{amscd}
\usepackage[all, cmtip]{xy}
\usepackage{mathtools}

\newtheorem{thm}{Theorem}[section]
\newtheorem{lem}[thm]{Lemma}
\newtheorem{prop}[thm]{Proposition}

\theoremstyle{defn}

\theoremstyle{remark}

\numberwithin{equation}{section}

\def \R  {\mathbb R}
\def \H {\mathcal H}
\def \F {\mathcal F}

\def \bt {\begin{tabular}{l}}
\def \et {\end{tabular}}

\title{First occurring singularities of functions on symplectic semi-space}
\date{}
\author{K. Kourliouros$^*$, M. Zhitomirskii\footnote{The work of both authors was supported by the Israel Science Foundation grant 510/12. \newline AMS Subject classification: 53D05, 58K50. \newline Key-words: symplectic semi-space, singularities, functional moduli.}
}

\begin{document}

\maketitle

\begin{abstract}
We explain how a classical theorem of Arnol'd and Melrose on non-singular functions on a symplectic manifold with boundary can be proved in few lines, and we use the same method to obtain a new result, which is a normal form with functional invariants for the first occurring singularities.
\end{abstract}

\section{Introduction and main results}
\label{sec-intr-results}

All objects in the paper are either $C^{\infty}$ or analytic germs at $0$.
By a symplectic semi-space  we mean the symplectic space $(\mathbb{R}^{2n},\omega )$ where $\omega $ is a symplectic form,  endowed with a smooth hypersurface  $\mathcal H$, which can be interpreted as the boundary.
The Darboux-Givental' theorem (see \cite{A3}) implies that all symplectic semi-spaces are locally equivalent.
We deal with the problem of classification of functions $f$ on a symplectic semi-space $(\R^{2n},\omega, \H)$.  Two function germs are called equivalent if they can be brought one to the other by a local symplectomorphism of $(\R^{2n}, \omega )$ which preserves the hypersurface $\H$.  An equivalent problem
is the classification of triples $(\omega,\H,f)$ with respect to the whole group of local diffeomorphisms.

\medskip

{\sc Notation}. In what follows $f$ and $h$ are function germs such that
\begin{equation}
\label{h}
\bt
$f(0)=0$, \ and \\
$h(0)=0, \ dh(0)\ne 0, \ \H =\{h=0\}$.
\et
\end{equation}

\noindent {\bf Definition}. {\it
The triple $(\omega , \H, f)$  is non-singular if $\{f,h\}(0)\ne 0$.}

\medskip

Here and in what follows $\{., .\}$ is the Poisson bracket with respect to
$\omega $.

\begin{thm}[Melrose \cite{M}, Arnol'd \cite{A1}]
\label{thm-Ar-Me}
Any non-singular triple $(\omega , \H, f)$ on $\R^{2n}$ is equivalent to
\begin{equation}
\label{nf-Ar-Me}
\bt
$\omega = \sum _{i=1}^ndx_i\wedge dy_i, \ \H=\{x_1=0\}, \ f=y_1$.
\et
\end{equation}
\end{thm}

In \cite{M} and \cite{A1} the proofs are given without details, but in fact
Theorem \ref{thm-Ar-Me} is a simple corollary of the Darboux
theorem on odd-symplectic (or quasi-symplectic) forms, and we find worth to give a very short proof in section \ref{sec-proofs}; the proofs of our much more difficult theorems on singular $(\omega , \H, f)$ are based on the same approach.

\medskip

The paper is devoted to first occurring singularities of $(\omega, \H, f)$ on $\R^{2n}$.

\medskip

\noindent {\bf Definition}.
{\it By $S_1$ we denote the singularity class consisting of $(\omega , \H, f)$ such that}
\begin{equation*}
\bt
for $n=1$: \ $\{f,h\}(0)=0$, \ $\{f,\{f,h\}\}(0)\ne 0$;
\\
for $n\ge 2$:  \ $\{f,h\}(0)=0$, \ $\{f,\{f,h\}\}(0)\ne 0$, \ $\{h,\{f,h\}\}(0)\ne 0$, \ $df\wedge dh(0)\ne 0$.
\et
\end{equation*}

It is easy to check that the choice of $h$ in (\ref{h}) in  the given  definitions is irrelevant: if they hold for $h$ then they also hold for $Qh$ where $Q$ is any non-vanishing function.
The difference in the definition of $S_1$ for $n=1$ and $n\ge 2$ is explained as follows: it is easy to see that for $n=1$ the assumptions
$\{f,h\}(0)=0$, \ $\{f,\{f,h\}\}(0)\ne 0$ imply that $\{h,\{f,h\}\}(0)\ne 0$ and $(df\wedge dh)(0)=0$.

\medskip

The problem of classification of triples $(\omega , \H, f)\in S_1$
was raised by R. B. Melrose in \cite{M} where
he studied a tied problem of classification of
triples $(\omega , \H,  \F)$, $\H=\{h=0\}, \ \F=\{f=0\}$
where $\H$ and $\F$ are two
smooth hypersurfaces and $\omega $ is a symplectic form.
Melrose studied the first occurring singularities of such triples,
called glancing hypersurfaces in a symplectic space,
which are distinguished by the same conditions as in the definition of $S_1$.
The main theorem in \cite{M} states that
in the $C^\infty $ category the glancing hypersurfaces in the symplectic space $\mathbb{R}^{2n+2}$ 
can be described by the normal form  
\begin{equation}
\label{Melrose-glancing-nf}
\bt
$\omega = dx\wedge dy + \sum _{i=1}^{n} dp_i\wedge dq_i$, \
$\H = \{y+x^2+p_1=0\}$, \ $\mathcal F=\{y=0\}$.
\et
\end{equation}
(without $dp_i\wedge dq_i$ in $\omega $ and without $p_1$ in $\H$ for $n=0)$. In \cite{M}, p. 176 he noticed that replacing the hypersurface $\F$ by a function $f$
makes the problem substantially different - starting from the first occurring singularities moduli occur.  Melrose claimed, and left to the reader to check, that in the family of triples
\begin{equation}
\label{Melrose-example}
\bt
$\omega = dx\wedge dy + \sum _{i=1}^{n}dp_i\wedge dq_i$, \
$\H = \{y+x^2+p_1=0\}$, \ $f = sy$
\et
\end{equation}
the parameter $s$ is an invariant.

\medskip

The problem of classifying singularities of $(\omega, \H, \F)$  is
tied with a number of other classification problems, many of them
were solved ( a good part - by V. I. Arnol'd and his school, c.f. \cite{A1}, \cite{A3}).
Nevertheless, the problem of classification of singularities of
functions in a symplectic semi-space remains, as we know, open.
The purpose of this paper is to fill in this gap.

\begin{thm}
\label{thm-n-1}
Any triple $(\omega , \H, f)\in S_1$ on $\R^2$ can be brought to
the normal form
\begin{equation}
\label{nf-n-1}
\bt
$(\omega , \H, f)_{g(y)}$:\\
$\omega = dx\wedge dg(y), \ \H=\{y+x^2=0\}, \ f = g(y)$,
\ \
$g^\prime (0)\ne 0.$
\et
\end{equation}
The triple (\ref{nf-n-1}) is equivalent to $(\omega , \H, f)_{\widetilde g(y)}$ if and only if \ $\widetilde g(y)=g(y)$.
\end{thm}

We have a functional invariant - a function of one variable $g(\cdot)$.
In \cite{K} the functional invariant was constructed in a canonical (coordinate-free) way using powerful tools from Gauss-Manin theory. In the same paper
another normal form
\begin{equation}
\label{nf-K}
\bt
$\omega = \phi(y+x^2)dx\wedge dy, \ \H=\{y=0\}, \ f=y+x^2$, \
$\phi(0)\ne 0$
\et
\end{equation}
with the functional invariant $\phi(\cdot )$ was obtained. Note that (\ref{nf-K})
can be easily obtained from (\ref{nf-n-1}) and the proof of
(\ref{nf-n-1}) in section \ref{sec-proofs} is elementary and takes just few lines.

\begin{thm}
\label{thm-n-from-2}
In the space of $(2n)$-jets of triples $(\omega , \H, f)\in S_1$
on $\R^{2n+2}$, $n\ge 1$, there exists an open set $U$ such that any triple $(\omega , \H, f)\in S_1$ with $j^{2n}(\omega , \H, f)\in U$
can be brought to the normal form
\begin{equation}
\label{nf-n-from-2}
\bt
$(\omega , \H, f)_{\mu , g(y), \phi(y,p,q)}$: \\
$\omega = dx\wedge dF(y,p,q) + \mu , \ \H=\{y+x^2=0\}, \ f = F(y,p,q)$
\\
$F(y,p,q)= g(y)+\sum_{i=1}^n(p_iy^{2i-2}+q_iy^{2i-1})+y^{2n}\phi(y,p,q)$, \ $g^\prime (0)\ne 0$,
\et
\end{equation}
where $\mu $ is a symplectic 2-form on $\R^{2n}(p,q)$. The  triple (\ref{nf-n-from-2})
is equivalent to
$(\omega , \H, f)_{\widetilde \mu , \widetilde g(y), \widetilde \phi(y,p,q)}$
if and only if $\widetilde \mu =\mu$, \ $\widetilde g(y)=g(y)$,\ $\widetilde \phi(y,p,q) = \phi(y,p,q)$.
\end{thm}

In both normal forms (\ref{nf-n-1}) and (\ref{nf-n-from-2}) the numerical invariant $g^\prime (0)$ has a simple canonical meaning
by the following lemma.

\begin{prop}
\label{lem-first-invariant}
In the classification of triples $(\omega , \H, f)\in S_1$
    the following number $\kappa $ is an invariant:
    \begin{equation*}
\bt
$\kappa =\frac{\{h,\{f,h\}\}(0)}{(\{f,\{f,h\}\}(0))^2}$
\et
\end{equation*}
In terms of normal forms (\ref{nf-n-1}) and (\ref{nf-n-from-2}) one has $\kappa = \frac{1}{2g^\prime (0)}$.
\end{prop}

\begin{proof}
Using that $h(0)=0$ and $\{f,h\}(0)=0$ it is easy to see that when multiplying $h$ by a non-vanishing function $Q$ one has
$\{f,\{f,h\}\}(0)\to Q(0)\{f,\{f,h\}\}(0)$ and
$\{h,\{f,h\}\}(0)\to Q^2(0)\{h,\{f,h\}\}(0)$. It follows that $\kappa $
does not depend of the choice of $h$. The formula $\kappa = \frac{1}{2g^\prime (0)}$
can be easily computed.
\end{proof}

\medskip

{\sc Remark}. Due to Lemma \ref{lem-first-invariant} we can do Melrose's homework:
the invariant $\kappa $ in his example (\ref{Melrose-example}) is equal to
$\frac{1}{2s}$, where $s=g'(0)$.

\medskip

Theorems \ref{thm-n-1} and \ref{thm-n-from-2} are proved, along with
much simpler and known Theorem \ref{thm-Ar-Me}, in the same way in section \ref{sec-proofs}. The proofs are very short provided that one uses the two lemmas in section \ref{sec-lemmas}, the proof of each one of them occupying only a few lines.
Probably the reader would ask why we do not use
the normal form $\omega = \sum dx_i\wedge dy_i, \ \H=\{x_1=0\}$  for
the pair $(\omega , \H)$ which seems a natural way to obtain normal forms for $(\omega , \H, f)$ by normalizing the function $f$ with respect
to local diffeomorphism preserving this normal form.
One can check that this way leads to very involved computations, and even if the computational obstacles can be resolved the proof would be very long.

\section{Auxiliary lemmas}
\label{sec-lemmas}

{\bf Definition}. {\it
A local odd-symplectic form on $\R^{2n+1}$ is the germ $\widehat w$  at $0$ of closed 2-form $\widehat \omega $ such that $\widehat \omega ^n(0)\ne 0$.}

\medskip

In the following lemma and its proof
\begin{equation}
\label{notations}
\bt
$p=(p_1,..., p_n), q = (q_1,..., q_n)$, $dp\wedge dq = \sum _{i=1}^ndp_i\wedge dq_i$.
\et
\end{equation}

\begin{lem}
\label{lem-cor-Darboux-odd-1}
  Any odd-symplectic form $\widehat \omega $ on $\R^{2n+1}(y,p,q)$ such that $(\widehat \omega ^n\wedge dy)(0)\ne 0$ can be brought to the form $dp\wedge dq$ by a local diffeomorphism preserving the coordinate $y$.
If $\Psi $ is a local diffeomorphism of  $\R^{2n+1}(y,p,q)$
which preserves $dp\wedge dq$ and the coordinate $y$ then $\Psi :(y,p,q)\to (y, A(p,q),B(p,q))$ where $(p,q)\to (A,B)$ is a local symplectomorphism of
$(\R^{2n}, dp\wedge dq)$.
\end{lem}

\begin{proof}
By the classical Darboux theorem for odd-symplectic forms (the proof can be found, for example, in \cite{Zh}) $\widehat \omega $ can be brought to
$dp\wedge dq$ by some local diffeomorphism $\Phi $. This diffeomorphism
brings $y$ to a function $Y(y,p,q)$; the assumption $(\widehat \omega ^n\wedge dy)(0)\ne 0$ implies $\frac{\partial Y}{\partial y}(0)\ne 0$.
 The local diffeomorphism $(y,p,q)\to (Y,p,q)$ brings $Y(y,p,q)$ to $y$ and preserves $dp\wedge dq$. It proves the first statement. To prove the second statement it suffices to note that any diffeomorphism which preserves
 both $y$ and the line field $\mathcal L$ generated by
$V=\frac{\partial }{\partial y}$ has the form given in the lemma, and that
$\mathcal L$ is invariantly related to $\widehat \omega = dp\wedge dq$:
$V$ is defined, up to multiplication by a function, by $V\rfloor \Omega = \widehat \omega ^n$ where $\Omega $ is a volume form.
 \end{proof}

\begin{lem}
\label{lem-Z-H-1}
Let $Z$ be a non-singular vector field on $\R^{2n+2}$ which has the simple tangency with
a smooth hypersurface $\H$. In suitable local coordinates $Z = \frac{\partial}{\partial x}, \ \H=\{y+x^2=0\}$. Any local diffeomorphism of $\R^{2n+2}$ which preserves this normal form preserves the coordinates $x$ and $y$.
\end{lem}

\begin{proof}
The given normal form is well-known, see for example
\cite{A0}. Let us prove the second statement. Any local diffeomorphism
preserving the vector field $\frac{\partial}{\partial x}$ changes $x$ to
$x+A(y,p,q)$ and $y$ to $B(y,p,q)$.  It preserves the hypersurface 
$y+x^2=0$ if and only if the function $(x+A)^2+B$ vanishes on the hypersurface $y=-x^2$. It means $x^2+B(-x^2,p,q)+A^2(-x^2,p,q)+ 2xA(-x^2,p,q)\equiv 0$.
Taking the even and the odd part with respect to $x$ we obtain $A=0, B=y$.
\end{proof}

\section{Proof of Theorems \ref{thm-Ar-Me}, \ref{thm-n-1}, \ref{thm-n-from-2}}
\label{sec-proofs}

In what follows we denote by $Z_f$ the Hamiltonian vector field
defined by $f$:
\begin{equation}
\label{Ham-vf}
Z_f \ \rfloor \ \omega = df.
\end{equation}
We use the coordinates $x,y$ in the 2-dim case. For higher dimensions
we work on $\R^{2n+2}(x,y,p,q)$ and we use notations (\ref{notations}).

\begin{lem}
\label{lem-omega-hat}
Let $Z_f = \frac{\partial }{\partial x}$. Then in the same coordinates
\begin{equation}
\label{eq-with-omega-hat}
\bt
\text{\rm  in the 2-dim case} $f = g(y)$, \ $\omega = dx\wedge dg(y)$
\\
\text{\rm for higher dimensions}:
\
$f = F(y,p,q)$, \ $\omega = dx\wedge dF(y,p,q)+\widehat \omega $
\\
\text{\rm where} $\widehat \omega $ \text{\rm is an odd-symplectic form on} $\R^{2n+1}(y,p,q)$.
\et
\end{equation}
\end{lem}

\begin{proof}
The given form of $f$ follows
from the equation $Z_f(f)=0$ which is a corollary of (\ref{Ham-vf}) with
$Z_f = \frac{\partial }{\partial x}$. This equation also implies
$\omega = dx\wedge dg(y)$ for $n=0$ and
$\omega = dx\wedge dF(y,p,q)+\widehat \omega $, \ $\frac{\partial }{\partial x}\rfloor \omega = 0$ for $n\ge 1$. Since $d\omega =0$ and
$\omega ^{n+1}(0)\ne 0$
it follows that $\widehat \omega $ is an odd-symplectic form.
\end{proof}

\smallskip

{\bf Proof of Theorem \ref{thm-Ar-Me}}. The condition $\{f,h\}(0)\ne 0$
means that $Z_f$ is transversal to $\H$.
Take coordinates in which $Z_f=\frac{\partial}{\partial x}, \ \H = \{x=0\}$. By Lemma \ref{lem-omega-hat}, $\omega $ and $f$ have form
(\ref{eq-with-omega-hat}). Changing the coordinate $y$ and using Lemma \ref{lem-cor-Darboux-odd-1} we obtain (\ref{nf-Ar-Me}), up to notations of the coordinates.

\medskip

{\bf Proof of Theorem \ref{thm-n-1}}.
The conditions $\{f,h\}(0)=0$, $\{f,\{f,h\}\}(0)\ne 0$  mean that $Z_f$ is a non-singular vector field which has the simple tangency with $\H$.
By Lemma \ref{lem-Z-H-1} the pair $(Z_f,\H)$ can be brought to the normal form $Z_f=\frac{\partial}{\partial x}, \ \H = \{y+x^2=0\}$.
By Lemma \ref{lem-omega-hat}, $\omega $ and $f$ have form
$\omega = dx\wedge dg(y), \ f = g(y)$ so that we have normal form (\ref{nf-n-1}) for the triple $(\omega , \H, f)$.  The fact that $g(y)$ is a functional invariant is a direct corollary of the second statement of Lemma \ref{lem-Z-H-1}.

\medskip

{\bf Proof of Theorem \ref{thm-n-from-2}}.
As in the proof of  Theorem \ref{thm-n-1}, we take coordinates in which
$Z_f=\frac{\partial}{\partial x}, \ \H = \{y+x^2=0\}$, and by
Lemma \ref{lem-omega-hat}, $\omega $ and $f$ have form
(\ref{eq-with-omega-hat}).
Now we will show that
the condition $\{h,\{f,h\}\}(0)\ne 0$ in the definition of the singularity class $S_1$ is equivalent to $(\widehat \omega ^n\wedge dy)(0)\ne 0$.
Let $h = y+x^2$.  Let $Z_h$ be the Hamiltonian vector field defined by $h$. The condition $\{h,\{f,h\}\}(0)\ne 0$ means $(Z_h(x))(0)\ne 0$. Indeed, since $Z_f=\frac{\partial }{\partial x}$ we have
$\{f,h\}=2x$ and $\{h, \{f,h\}\} = 2Z_h(x)$. 
To find $Z_h(x)$ we have by definition
$Z_h \ \rfloor \ (dx\wedge dF(y,p,q) + \widehat \omega ) = dh = dy + 2xdx$.
Express this equation in the form
\begin{equation}
\label{eq-t}
\bt
$Z_h(x)\cdot dF(y,p,q) = dy - Z_h\rfloor \widehat \omega + (2x+(Z_h\ \rfloor\  dF(y,p,q))dx$
\et
\end{equation}
Since $\widehat \omega ^{n+1}=0$ we have $(Z_h\rfloor \widehat \omega )\wedge \widehat \omega ^n=0$. Therefore taking the external product of
(\ref{eq-t}) with $\widehat \omega ^n\wedge dx$ we obtain
$Z_h(x)\cdot dF(y,p,q)\wedge \widehat \omega ^n\wedge dx  = dy\wedge \widehat \omega ^n\wedge dx$ and it follows
$dy\wedge \widehat \omega ^n =
Z_h(x)\cdot (dF(y,p,q)\wedge \widehat \omega ^n)$. Since $\omega ^{n+1}\ne 0$ we have $dF(y,p,q)\wedge \widehat \omega ^n\ne 0$. Therefore
the condition $(dy\wedge \widehat \omega ^n)(0)\ne 0$ is the same as
$(Z_h(x))(0)\ne 0$ which, as is shown above, is the same as $\{h,\{f,h\}\}(0)\ne 0$. 
\smallskip

Now we can use Lemma \ref{lem-cor-Darboux-odd-1} which allows to bring $\widehat \omega $ to
$dp\wedge dq$ preserving $Z_f = \frac{\partial }{\partial x}$ and
$\H=\{y+x^2=0\}$. We obtain the normal form
\begin{equation}
\label{nf-pre-final}
\bt
$(\omega , \H, f)_{F(y,p,q)}: $
\\
$\omega = dx\wedge dF(y,p,q) + dp\wedge dq, \ \H=\{y+x^2=0\}, \ f = F(y,p,q)$.
\et
\end{equation}
The next step is using the second statements of Lemmas \ref{lem-cor-Darboux-odd-1} and \ref{lem-Z-H-1}. They imply that $(\omega , \H, f)_{F(y,p,q)}$ is equivalent to $(\omega , \H, f)_{\widetilde F(y,p,q)}$ if and only if
the functions $F(y,p,q)$ and $\widetilde F(y,p,q)$ can be brought one to the other by a symplectomorphisms of $\R^{2n}(p,q)$.
Therefore reducing (\ref{nf-pre-final}) to exact normal form is the same problem as classification of functions $F(y,p,q)$ with respect to  symplectomorphisms of $\R^{2n}(p,q)$.
An equivalent problem is as follows:

\medskip

 \noindent (*) classification of pairs $(\mu, F(y,p,q))$ where $\mu $ is any symplectic form on $\R^{2n}(p,q)$ with respect to all diffeomorphisms of the form $(y,p,q)\to (y,\Phi (p,q), \Psi (p,q))$.
%

\medskip

Let $F(y,p,q)=g(y) +F_0(p,q) + F_1(p,q)y + \cdots $ where $F_i(0,0)=0$. Assuming that the functions $F_0,..., F_{2n-1}$ are differentially independent at $0$, which defines the open set $U$ in Theorem \ref{thm-n-from-2}, we obtain, in the problem (*), the exact
normal form with
$F(y,p,q)=g(y)+p_1+q_1y + p_2y^2+q_2y^3\cdots + q_ny^{2n-1}+ \phi(y,p,q)y^{2n}$  where $g(y)$ and $\phi(y,p,q)$ are functional invariants, and the symplectic form $\mu $ in the pair $(\mu, F(y,p,q))$ is a functional invariant too.
This exact normal form leads to the exact normal form (\ref{nf-n-from-2}).

\end{document}